\documentclass[pdflatex,sn-mathphys-ay]{sn-jnl}

\usepackage{graphicx}%
\usepackage{multirow}%
\usepackage{amsmath,amssymb,amsfonts}%
\usepackage{amsthm}%
\usepackage{mathrsfs}%
\usepackage[title]{appendix}%
\usepackage{xcolor}%
\usepackage{textcomp}%
\usepackage{manyfoot}%
\usepackage{booktabs}%
\usepackage{algorithm}%
\usepackage{algorithmicx}%
\usepackage{algpseudocode}%
\usepackage{listings}%
\usepackage{natbib}
\usepackage{comment}
\usepackage{eurosym}

\theoremstyle{thmstyleone}%
\theoremstyle{thmstyletwo}%

\theoremstyle{thmstylethree}%

\DeclareMathOperator*{\argmin}{arg\,min}

\begin{document}

\title[Article Title]{A QUBO-Based Optimization Framework for ATM Cash Replenishment Scheduling}


\author*[1]{\fnm{Riccardo} \sur{Aiolfi}}\email{r.aiolfi@reply.it}

\author[1]{\fnm{Giulia} \sur{Montani}}

\author*[2]{\fnm{Valeria} \sur{Zaffaroni}}\email{valeria.zaffaroni@intesasanpaolo.com}

\author[1]{\fnm{Luca} \sur{Toffanetti}}\email{l.toffanetti@reply.it}



\author[2]{\fnm{Davide} \sur{Veronelli}}\email{davide.veronelli@intesasanpaolo.com}

\author[2]{\fnm{Francesca} \sur{Cibrario}}\email{francesca.cibrario@intesasanpaolo.com}

\author[1]{\fnm{Davide} \sur{Caputo}}\email{da.caputo@reply.it}

\author[2]{\fnm{Davide} \sur{Corbelletto}}\email{davide.corbelletto@intesasanpaolo.com}

\affil[1]{\orgdiv{Data Reply s.r.l.}, \orgaddress{\street{Corso Francia 110}, \city{Torino}, \postcode{10143}, \state{Italy}}}

\affil[2]{\orgdiv{Intesa Sanpaolo S.p.A.},  \orgaddress{\street{Piazza San Carlo 156}, \city{Torino}, \postcode{10121}, \state{Italy}}}


\abstract{
The management of cash replenishment in Automated Teller Machine (ATM) networks requires scheduling recharges in order to minimize operational costs while maintaining high service levels and avoiding cash-outs, under uncertain and time-varying withdrawal demand.
This work formulates the ATM cash replenishment problem through a \emph{Quadratic Unconstrained Binary Optimization} (QUBO) model, which naturally captures nonlinear cost interactions, while incorporating operational constraints through penalty terms.
The objective function combines fixed and variable replenishment costs with co-location discounts, as well as penalties for a late replenishment that could cause a service interruption.
The resulting QUBO instances are solved using MegaQUBO, a GPU-accelerated QUBO solver.
An empirical evaluation on a real dataset of 276 ATMs located in Italy, covering four representative months of 2022 (April, May, October, and November), benchmarks the proposed approach against a threshold-based operational policy.
Results show consistent cost reductions of approximately 15\%--18\% while maintaining an excellent average service level (around 99.8\%--99.9).
Overall, the study demonstrates that QUBO-based optimization, coupled with GPU-based solving, can provide a practically deployable decision-support tool for large-scale ATM cash logistics.}




\maketitle

\section{Introduction}\label{sec1}
The management of cash within logistics networks and the optimization of cash flow present significant operational challenges for financial institutions worldwide. For networks of Automated Teller Machines (ATMs), the fundamental problem is to schedule cash replenishments in a way that minimizes total operational costs while maintaining high service levels and avoiding cash-outs that negatively impact customer satisfaction. The complexity of this problem stems from multiple interacting factors: uncertain and time-varying withdrawal demand, the combinatorial nature of routing decisions, capacity constraints on replenishment vehicles, and non-linear cost structures when servicing multiple machines in proximity.

Classical approaches to ATM cash replenishment have traditionally relied on simple heuristics or adaptations of well-known operations research methodologies. The most common practical approach employs threshold-based policies, where replenishment is triggered when a machine's projected cash balance falls below a predetermined level. Although straightforward to implement, this strategy fails to optimize costs across the network and does not account for the complex interdependencies between replenishment decisions. More sophisticated methods have adapted techniques from inventory control theory and vehicle routing optimization, formulating the problem as variants of the Multi-Depot Vehicle Routing Problem with Time Windows (MDVRPTW). Moreover, Mixed Integer Linear Programming (MILP) formulations have been successfully applied to small and medium-scale instances, while metaheuristic approaches—including genetic algorithms, simulated annealing, and tabu search—have provided near-optimal solutions for larger networks where exact methods become computationally prohibitive. Despite these advances, the ATM cash replenishment problem remains challenging due to its intrinsic computational complexity. Real-world instances often involve non-linear cost interactions that are difficult to capture in traditional linear formulations. For example, servicing multiple co-located ATMs in a single visit may yield significant cost discounts, creating quadratic relationships between decision variables that standard MILP solvers handle inefficiently.

This paper introduces a novel approach by formulating the ATM cash replenishment problem as a Quadratic Unconstrained Binary Optimization (QUBO) problem. 
The QUBO framework, indeed, naturally accommodates the quadratic cost interactions inherent in the problem, i.e. discounts for joint service, and allows complex operational constraints — such as vehicle capacity limits, time windows, and service level requirements — to be incorporated through penalty terms in the objective function. Moreover, QUBO is the native formulation for quantum annealing devices and hybrid quantum-classical solvers, which have shown promise in solving large-scale combinatorial optimization problems that are intractable for conventional methods. By recasting the ATM replenishment problem in the QUBO framework, we create a pathway to leverage emerging quantum computing technologies and advanced classical optimization techniques specifically designed for quadratic binary problems. In this work, we exploit this formulation by solving the resulting QUBO instances with MegaQUBO \citep{megaqubo}, a Data Reply's proprietary GPU-accelerated QUBO solver, which has already been employed to enable large-scale solutions, for example in \cite{Logistic} and \cite{RiskModel}.
Furthermore, this approach has been deployed in a production environment, thereby supporting the validity of the reported results.

The remainder of this work is structured as follows. Section~\ref{sec2} reviews the relevant literature on ATM cash management, classical optimization approaches, and recent applications of quantum computing and QUBO formulations in logistics. Section~\ref{sec3} presents the problem formulation and the proposed QUBO model for ATM cash replenishment. Section~\ref{sec4} describes the solution approach, the experimental setup and presents the results, comparing the proposed approach with the threshold-based operational strategy. Finally, Section~\ref{sec5} concludes with a discussion of results, practical implications, and directions for future research.
\section{Related works}\label{sec2}
The optimization of ATM cash replenishment has been studied extensively through various methodological lenses, spanning classical optimization approaches and stochastic methods.
A first stream of work focuses on exact optimization formulations for replenishment planning. ATM cash replenishment was studied in \citep{ozer2020} by comparing integer linear programming and dynamic programming approaches, emphasizing the computational burden that arises as the network scale increases and showing that dynamic programming can offer a more efficient solution procedure for the considered setting. A complementary line of research highlights the role of forecast-informed decision making: for example, in \citep{cruzreyes2015}, demand-forecast information was integrated into the replenishment optimization process, illustrating how forecast-driven planning can influence replenishment choices and overall system performance.

In parallel, quantum computing has motivated new modeling paradigms for large-scale combinatorial optimization. Within this context, QUBO formalism has emerged as a convenient representation due to its ability to encode binary decisions and quadratic interactions directly in the objective function. QUBO models have been explored on quantum annealing hardware for dynamic routing problems \citep{hussain2020} and extended to quadratic and higher-order unconstrained formulations for complex rescheduling tasks \citep{zaman2022}. Closer to the financial logistics domain, in \citep{diezvalle2023}, a multi-objective variational quantum optimization for cash management was investigated, suggesting the potential of quantum optimization methods when multiple objectives and constraints must be handled simultaneously.

While classical optimization methods have been extensively studied for ATM cash replenishment and quantum computing has shown promise in various logistics domains, the application of QUBO formulations specifically to ATM cash replenishment remains largely unexplored. This gap is significant because ATM replenishment exhibits structural properties that align naturally with the QUBO framework: binary decision variables (service/no-service), quadratic cost interactions (joint servicing discounts), and the possibility to encode complex constraints through penalty terms. To bridge this gap and establish a rigorous foundation for our approach, Section~\ref{sec3} introduces the QUBO framework and presents our problem-specific formulation for cash flow optimization in ATM logistics.
\section{Methodological Background}\label{sec3}
This section surveys the theoretical foundations and practical formulation of the optimization approach developed for ATM cash replenishment scheduling. While the QUBO formulation constitutes the principal innovative contribution of this work, its role and significance can only be fully appreciated when placed within the broader end-to-end pipeline of which it forms the decision-making core.  

\subsection{End-to-end Pipeline}
The pipeline operates on a daily cycle and begins with the ingestion of operational data required by downstream components such as real-time cash residual of each ATM, machine-specific parameters (maximum capacity and geographic co-location groupings) and historical transaction records spanning a sufficient window for demand estimation. The ingested data feed a time-series forecasting module that produces estimates of expected daily withdrawal volumes for each ATM over a planning horizon. At the center of the pipeline lies the QUBO-based optimization model, which determines the optimal replenishment schedule across all ATMs and all days in the planning window. A binary decision variable is assigned to each (ATM, day) pair: a value of one indicates that a replenishment intervention is to be executed for the corresponding machine on the corresponding day, while a value of zero indicates that no intervention is planned. Since the replenishment service ticket is generated whenever a machine's observed cash residual falls below a prescribed value, the replenishment schedule suggested by the QUBO is translated into a set of daily activation thresholds. The translation logic distinguishes two cases: if the QUBO solution schedules a replenishment for a given ATM on day $t+1$, then the activation threshold for that ATM on day $t$ is set equal to the machine's maximum capacity. Because the real-time residual cannot exceed this capacity, the threshold condition is guaranteed to be satisfied, and a service ticket is opened on day $t$ in preparation for execution on the following day. Conversely, if no replenishment is scheduled for a given ATM on the next day, the threshold is computed as the cumulative sum of forecasted daily withdrawals over the remaining days of the horizon until a replenishment intervention could be scheduled ---including weekends and non-working days during which interventions cannot be performed---augmented by a safety contingency margin that absorbs demand variability and reduces the probability of unplanned cash-outs. The operational system evaluates this condition daily: the real-time cash residual of each ATM is compared against its corresponding threshold, and if the residual falls below this value, a replenishment ticket is triggered in the service management platform.  

Having established the architectural context within which the QUBO formulation operates, in the following subsections we formalize the QUBO framework and discuss its expressiveness for modeling combinatorial structures, including hard constraints via penalty encodings. We then introduce a QUBO formulation tailored to ATM cash flow optimization, detailing the construction of the objective function and the incorporation of feasibility and service level requirements as penalty terms. This background sets the stage for the subsequent empirical evaluation, in which the proposed QUBO model is solved and benchmarked against the threshold-based operational strategy. 
\subsection{Quadratic Unconstrained Binary Optimization}\label{subsec31}
QUBO is a canonical formulation for discrete optimization problems with binary decision variables. Given $n$ Boolean variables collected in a vector $\mathbf{x} \in \{0,1\}^n$ and a $n \times n$ upper triangular real-valued matrix $Q$, a QUBO instance can be written as
\begin{equation}
\hat{\mathbf{x}} = \argmin_{\mathbf{x} \in \{0,1\}^n} \mathbf{x}^\top Q \mathbf{x}
= \argmin_{\mathbf{x} \in \{0,1\}^n} \sum_{i=0}^{n-1} \sum_{j=i}^{n-1} x_i \, Q_{ij} \, x_j 
\end{equation}
Despite its compact form, the QUBO framework is expressive enough to capture a broad set of NP-hard problems through suitable encodings of the objective and constraints.

Hard constraints are commonly incorporated in QUBO by augmenting the objective with penalty terms. Under this approach, a constrained optimization problem is recast as an unconstrained quadratic objective of the form
\begin{equation}
\min_{\mathbf{x}\in\{0,1\}^n}\; f(\mathbf{x}) \;+\; \sum_{k} \lambda_k\,\pi_k(\mathbf{x})
\end{equation}
where $f(\mathbf{x})$ is the primary cost function and each $\pi_k(\mathbf{x})\ge 0$ measures the violation of the constraint $k$, attaining its minimum (typically zero) on feasible assignments. The penalty weights $\lambda_k>0$ are chosen sufficiently large so that any constraint-violating configuration is dominated by at least one feasible configuration in terms of objective value. This penalty-based reformulation is widely used in practice when constructing QUBO (or, more generally, binary quadratic) models; see, for example, D-Wave's discussion of penalty models and problem reformulation \citep{dwave_penalty,dwave_reformulating}. 

A frequently used building block is the at-most-one (AMO) constraint, which enforces that from a set of binary indicators $\{x_{1},\ldots,x_{T}\}$ at most one can be selected:
\begin{equation}
\sum_{t=1}^{T} x_t \le 1
\end{equation}
In QUBO, AMO can be encoded via a quadratic penalty that assigns a positive cost to every pair of simultaneously active variables \citep{Alidaee2008},
\begin{equation}
P \cdot \sum_{t_1=1}^{T} \sum_{t_2=t_1+1}^{T} x_{t_1}x_{t_2}
\end{equation}
with $P>0$, a constant sufficiently large to be properly tuned. 

\subsection{QUBO Formulation for Cash Flow Optimization}\label{subsec32}
We consider the problem of scheduling cash replenishments for $N$ machines (ATMs) over a planning horizon of $T$ days, with the objective of minimizing total operational costs while maintaining a high service level. The daily replenishment cost for each machine comprises a fixed component and a variable component, where the latter increases with higher cash residuals on the machine. More specifically, in the considered setting, replenishment operations are carried out by fully replacing the residual cash with an amount corresponding to the maximum capacity of the associated ATM.
While fixed costs are related to the number of replenishment, variable costs are determined by two components: the volume of cash loaded into the ATM, which is weighted by a multiplicative coefficient, and the residual cash that is withdrawn and subsequently counted, which is subject to a distinct multiplicative coefficient. Additionally, when more machines located at the same address are replenished on the same day, a discount is applied. This discount introduces a nonlinearity that can be encoded within a Quadratic Unconstrained Binary Optimization (QUBO) framework. Considering $U$ the set of all ATMs that must be recharged and $A$ the set of their distinct addresses, the general cost function could be defined as:
\begin{equation}
\sum_{i\in U} \left( F + c \cdot C_i + s \cdot S_i \right) - \sum_{a\in A} (N_a - 1) \cdot D
\end{equation}
\noindent where the positive parameters are defined as follows:
\begin{itemize}
    \item $F$: fixed cost associated with the replenishment of a single ATM; 
    \item $C_i$: capacity of cash of ATM $i$;
    \item $c$: unit cost coefficient associated with the preparation and loading of cash to be used for replenishment;
    \item $S_i$: residual cash in ATM $i$ at the time of replenishment;
    \item $s$: unit cost associated with handling and counting the residual cash in the cash processing facility;
    \item $N_a$: number of distinct ATMs replenished at address $a$;
    \item $D$: constant cost saving obtained when multiple ATMs at the same address are replenished simultaneously.
\end{itemize}
The first summation captures the gross replenishment cost across all serviced machines, while the second summation implements the co-location discount: for each address at which more than one ATM is serviced simultaneously, a reduction proportional to $(N_a - 1)$ is applied, reflecting the shared logistical overhead. This expression represents the true operational cost structure as defined by the contractual service agreement. The subsequent QUBO formulation approximates this objective within the constraints imposed by the quadratic binary framework, as detailed in the following subsections.

To model this problem, we define an objective function $H(\mathbf{x})$ which is the sum of the operational cost function, discount terms for co-location, penalty terms for feasibility and service level constraints:
\begin{equation} 
H(\mathbf{x}) = \lambda_1H_{\text{cost}}(\mathbf{x}) + \lambda_2H_{\text{discount}}(\mathbf{x}) + \lambda_3H_{\text{feasibility}}(\mathbf{x}) + \lambda_4H_{\text{service}}(\mathbf{x}) 
\end{equation}
where 
\begin{equation}\label{eq:xbinary} 
\mathbf{x} = (x_{11}, ... ,x_{1T}, x_{21}, ... , x_{it}, ... , x_{NT}) 
\end{equation}
is a vector of binary decision variables indicating whether machine $i$ is replenished on day $t$ and $\{\lambda_{k}\}_{k=1,...,4}$ are strictly positive penalty weights that need to be tuned appropriately. This formulation allows for the incorporation of nonlinear cost interactions, address discounts, and service level requirements within a unified binary optimization framework, suitable for quantum-inspired solvers.

\subsection{Objective function}\label{subsec33}
In order to reformulate this scheduling problem within a QUBO framework, two simplifying assumptions are introduced. 

First, since the exact time at which a replenishment is physically executed during the day---i.e., the moment at which the service provider delivers the cash to the machine---is not known in advance, we conservatively assume that all replenishment operations take place at the end of the day. This guarantees that the ATM must be able to sustain demand throughout the entire day on which the replenishment is scheduled, thereby addressing the worst-case scenario and avoiding any dependence on uncertain delivery timing.

Secondly, the residual cash depends recursively on previous withdrawals and on whether a replenishment has occurred on any earlier day within the horizon. Modeling this recursion explicitly within a QUBO formulation would introduce higher-order interactions among decision variables across multiple time steps, breaking the quadratic structure. To retain QUBO compatibility and since the pipeline is executed every day, we impose that each ATM can be replenished \emph{at most once} within the planning horizon. Under this constraint, the impact of scheduling a replenishment on a given day $t$ can be evaluated by precomputing the cash trajectory that would materialize if no replenishment were performed at all. We therefore introduce the projected end-of-day balance under the \textit{no-replenishment scenario}:
\begin{equation}
\tilde{S}_{i,t} = \max(0,R_{i,1} - \sum_{j=1}^{t} W_{i,j})\,
\end{equation}
where $R_{i,1}$ denotes the known residual at the beginning of the first day in the planning horizon for the ATM $i$ and $W_{i,j}$ represents the forecasted withdrawal on day $j$ for the ATM $i$. 

After having introduced these premises, we can formalize the cost and discount components of the objective function for our optimization model. 
The total recharge cost for $N$ machines over a time horizon $T$ is
\begin{equation}\lambda_1\cdot\underbrace{\sum_{t=1}^T 
\sum_{i=1}^{N} x_{i,t} \, \big( F + c \cdot C_{i} + s \cdot \tilde{S}_{i,t} \big)}_{H_{\text{cost}}(\mathbf{x})}
+ \lambda_2\cdot \underbrace{\big(- D\big) \cdot \sum_{t=1}^T \sum_{i=1}^{N} \sum_{\substack{j=1 \\ j > i}}^{N} 
v_{i,j} \, x_{i,t} \, x_{j,t}}_{H_{\text{discount}}(\mathbf{x})}
\end{equation}
where $x_{i,t}$ is the binary decision variable 
and $v_{i,j}$ a fixed entry equal to $1$ if machine $i$ and $j$ are at the same address, $0$ otherwise. The discount term appearing on the right-hand side of the last expression is a mathematical implementation of the co-location discount applied when multiple ATMs at the same physical address are replenished on the same day. This formulation has been specifically adapted to comply with the quadratic structure required by the QUBO framework, in which the objective function must be expressible as a polynomial of degree at most two in the binary decision variables. This term does not always coincide exactly with the actual contractual discount; however, it constitutes a faithful approximation that captures the essential economic rationale of the co-location incentive—namely, encouraging the simultaneous servicing of multiple machines at the same site in order to reduce the total number of replenishment operations and the associated logistical costs. In particular, when exactly two ATMs located at the same address are scheduled for replenishment on the same day, this implementation corresponds precisely to the real discount. When the number of co-located ATMs replenished simultaneously exceeds two, the implemented discount slightly overestimates (in absolute value) the actual contractual reduction. Nevertheless, the cost function parameters can be straightforwardly recalibrated should the optimizer produce solutions that appear imbalanced or operationally inadequate.

\subsection{Constraints}\label{subsec34}
To ensure that the optimization model produces practical and reliable replenishment schedules, two categories of constraints are introduced: feasibility constraints are imposed to align the model with operational realities, while service level constraints are formulated to guarantee a high degree of machine availability, preventing cash-out events. In the QUBO framework, these constraints are incorporated as penalty terms, guiding the solver toward solutions that are both operationally valid and meet service quality standards.

\subsubsection{Feasibility Constraint}\label{subsubsec341}
The primary feasibility constraint is the single replenishment constraint. Historical data indicates that it is rare for a single machine to be replenished more than once within such a timeframe. Therefore, to simplify the model and eliminate complex recursive dependencies in cash balance calculations, we enforce that each ATM can be replenished at most once during the optimization horizon ($\sum_{t=1}^T x_{i,t} \le 1 \,\,\forall i\in\{1,...,N\}$). This requirement corresponds to an at-most-one (AMO) constraint and is implemented via a quadratic penalty that assigns a positive cost to any pair of selected replenishment days for the same ATM:
\begin{equation}
     \lambda_3 \cdot \underbrace{\sum_{i=1}^N \sum_{t_1=1}^T \sum_{t_2 = t_1 + 1}^T x_{i,t_1} x_{i,t_2}}_{H_{\text{feasibility}}(\mathbf{x})} \,
\end{equation}
The corresponding penalty weight must be chosen in such a way as to be significantly greater than any possible cost saving, ensuring that any solution violating this constraint is energetically unfavorable. This term is zero if an ATM is replenished zero or one time, and positive otherwise.

\subsubsection{Service Level Constraint}\label{subsubsec342}
The service level requirement aims to prevent cash-out events by ensuring that the expected cash balance of each ATM remains above a contingency level $\text{Fix}_i$ throughout the planning horizon. We define the cumulative shortage penalty associated with replenishing on day $t$ as 

\begin{equation}
\Gamma_{i,t} = \min\{0,\; \tilde{S}_{i,t} - \text{Fix}_i\}
\end{equation} 

\noindent which is a precomputable constant for each pair $(i,t)$. Moreover, we introduce a coefficient $p_{i,t} > 0$ that can be tuned to encourage recharging at the earliest time the ATM might become empty.
The service level term is then encoded as a linear penalty in the QUBO objective:
\begin{equation}
 \lambda_4 \cdot \underbrace{\sum_{i=1}^{N}\sum_{t=1}^{T} p_{i,t}
 \Gamma_{i,t}\, x_{i,t}}_{H_{\text{service}}(\mathbf{x})} \,
\end{equation}
Given that the structural constraint permits at most one replenishment event per ATM within the planning horizon, the penalty mechanism must selectively favor the specific $(i, t)$ pair corresponding to the first day on which a cash shortage would occur. To achieve this, the  coefficient $p_{i,t}$ must be sufficiently large on that first critical day so as to render it energetically more favorable than any alternative replenishment day. Indeed, for a fixed ATM, $\Gamma_{i,t}$ yields a sequence of zero values for all days preceding the projected shortage  and strictly negative values thereafter. The role of $p_{i,t}$  is therefore to ensure that the first non-zero entry in this sequence dominates all subsequent ones, thereby steering the optimizer toward scheduling the replenishment intervention on the earliest day at which service continuity would otherwise be compromised.

\section{Results and Discussion}\label{sec4}

\subsection{Dataset Specifications and Preprocessing}\label{subsec41}
The empirical basis for this study is a dataset detailing the operations of 276 ATMs distributed across Italy. The data spans four representative months of activity from 2022: April, May, October, and November. This dataset includes historical cash withdrawal information and machine-specific attributes such as location and operational constraints, i.e., ATM with or without Saturday Service. Since actual future withdrawal amounts are not available at the time of optimization, we introduce a time-series forecasting model that leverages historical data to produce estimates of expected daily withdrawals over the planning horizon.

A significant challenge in preparing the data for the QUBO model was standardizing the replenishment schedule across all ATMs, given their heterogeneous operational constraints. Although replenishment service tickets can only be issued from Monday to Friday, the actual service day depends on the ATM's specific service agreement. If a ticket is created on a Friday, an ATM with Saturday service availability will be replenished the next day. In contrast, an ATM without Saturday service will be replenished on the following Monday.

To maintain a uniform structure for the QUBO model, we adopted a fixed number of time variables for the planning horizon, but the meaning of these variables shifts depending on the ATM type: 
\begin{itemize} 
    \item \textbf{For ATMs with Saturday Service:} the time variables $t$ in a given week correspond to the set {Tuesday, Wednesday, Thursday, Friday, Saturday}; 
    \item \textbf{For ATMs without Saturday Service:} the time variables $t$ correspond to the set {Tuesday, Wednesday, Thursday, Friday, Monday}. 
\end{itemize}

This approach results in ``dual-meaning'' time indices. For example, within a 9-day rolling horizon, a specific time index might represent a Saturday for one type of ATM but a Monday for another. This ambiguity has two critical implications that must be handled during data preprocessing.
\begin{enumerate}
    \item\textbf{Withdrawal Forecast Adjustment.} To prevent cash-out events over non-service periods, withdrawal forecasts must be aggregated. For an ATM without Saturday service, the forecasted withdrawals for Saturday and Sunday are added to Friday's forecast. This ensures that a replenishment scheduled for Friday leaves the machine with sufficient cash to operate throughout the weekend until the next possible service on Monday.
    \item \textbf{Co-location Discount Constraint.} The cost-saving discount for replenishing multiple ATMs at the same address on the same day is affected. This discount can only be applied if the machines are serviced on the same actual calendar day. For a dual-meaning time index, the discount is only valid if the indexes of the co-located ATMs represent the same day. For example, if the index is Saturday for one ATM and Monday for another, the discount cannot be applied.
\end{enumerate}

Crucially, these necessary adjustments are implemented entirely within the pre-processing stage. This strategy allows the core QUBO model to remain structurally consistent and uniform for all ATMs, without requiring any modifications to handle the operational complexities arising from different service schedules.

\begin{figure}[h]
\centering
\includegraphics[width=\textwidth]{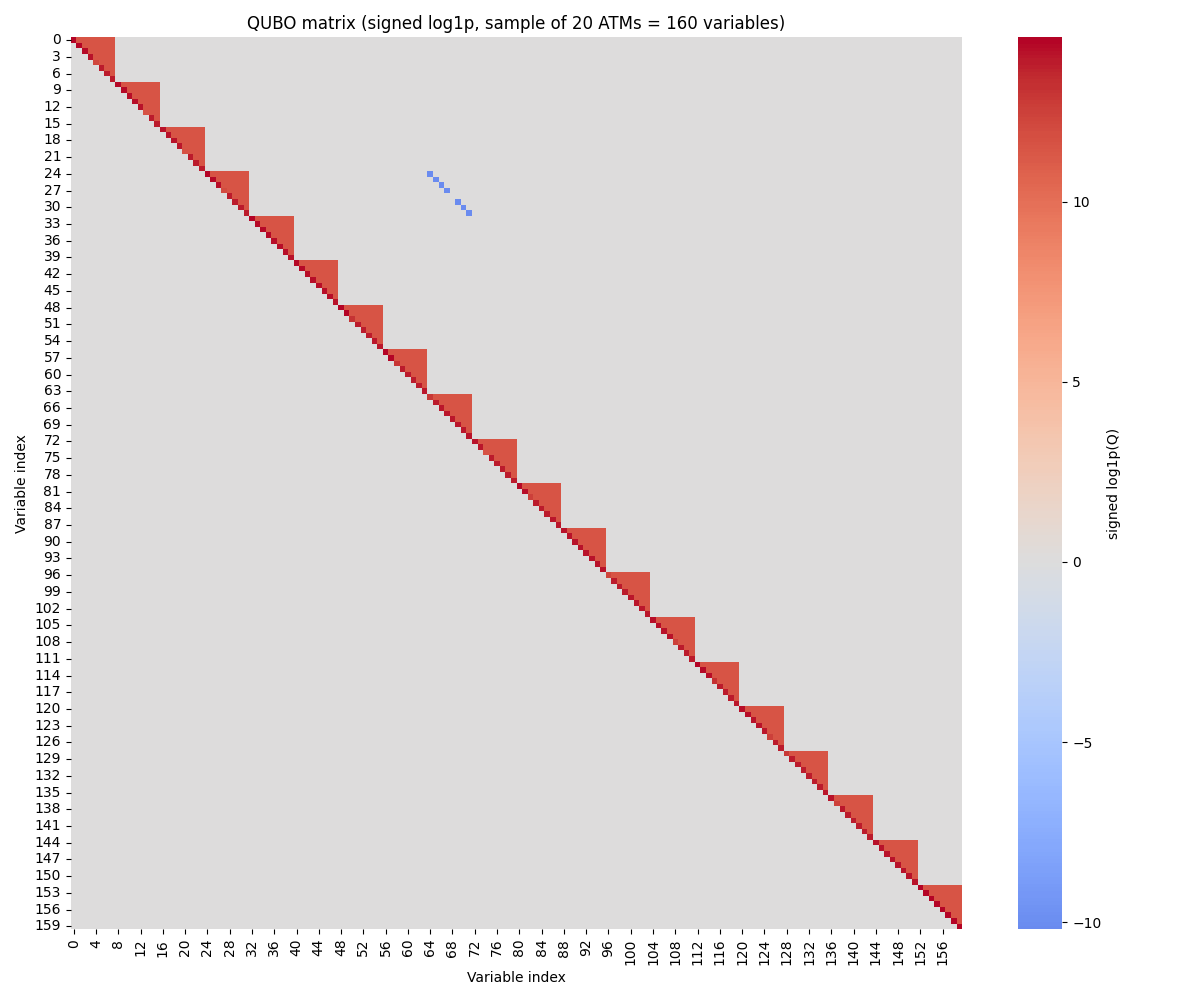}
\caption{Snapshot of the QUBO matrix $Q$ for 20 ATMs over an 8-day horizon (160 variables) with one co-location interaction pair between ATM 4 and 9. Values are displayed using a signed $\log(1+|Q_{ij}|)$ scaling to emphasize both sparsity and coefficient magnitudes.}
\label{fig:qubo_matrix}
\end{figure}

Fig.~\ref{fig:qubo_matrix} reports a snapshot of the rescaled QUBO coefficient matrix $Q$ for a reduced example consisting of 20 ATMs over an 8-day planning horizon (160 binary variables) with a single co-location interaction pair between ATM 4 and 9. As presented in section~\ref{subsec32}, the overall objective function is expressed as $H(\mathbf{x})=\mathbf{x}^\top Q\mathbf{x}$, where each variable corresponds to a replenishment decision $x_{i,t}$. With the ordering that groups the $T$ time variables of each ATM consecutively as showed in \eqref{eq:xbinary}, the matrix exhibits a characteristic block-diagonal structure. Linear terms from $H_{\text{cost}}(\mathbf{x})$ and $H_{\text{service}}(\mathbf{x})$ contribute to diagonal entries, whereas the AMO feasibility penalty $H_{\text{feasibility}}(\mathbf{x})$ induces dense within-ATM off-diagonal couplings across different days, forming the visible light red ATM-specific blocks. Cross-ATM couplings are sparse, negative (in blue) and arise only from the co-location discount term, which links variables $x_{i,t}$ and $x_{j,t}$ when $v_{i,j}=1$ and the same day $t$ is selected. 
In this example, we consider ATM 4 with Saturday Service and ATM 9 without Saturday Service. If the starting day is Tuesday, the discount is not applied on the 5th day for the "dual-meaning" time indices: Saturday for ATM 4 and Monday for ATM 9.
					
\subsection{Simulation Setup and Benchmark Solution}\label{subsec42}
To quantitatively assess the effectiveness of the proposed QUBO model, its performance is benchmarked against the threshold-based operational strategy within a controlled simulation environment. This approach ensures a direct comparison by removing external sources of noise, such as delays on recharges. The threshold-based replenishment strategy, serving as the baseline for our evaluation, operates on a simple, reactive threshold-based logic. In this model, a replenishment for an ATM is automatically triggered as soon as its residual cash balance falls below a threshold based only on the forecasted withdrawals. 
While straightforward and robust baseline, this method does not manage the overall cost optimization over a time horizon and does not take care of possible discounts due to recharging ATMs at the same address.
The comparison is performed using a simulation based on four months of real-world data from 2022, allowing for a realistic yet controlled evaluation. All cost parameters are expressed in euros. In the experiments of this study we set the fixed replenishment cost to $F=52$, the variable cost coefficient to $c=0.000069$, the residual-cash processing coefficient to $s=0.000333$, and the co-location discount to $D=26.69$. These parameter values are selected to ensure consistency with realistic operational settings and empirical cost structures.

\subsection{Evaluation Criteria}\label{subsec43}
The benchmarking methodology relies on a retrospective evaluation framework in which the actual realized withdrawal sequence—available only ex post—is used to assess the performance of each strategy. Specifically, given a historical period for which real data are known, we simulate the execution of both the proposed QUBO-based scheduling strategy and the  threshold-based policy, and we reconstruct the resulting cash residual trajectories for each ATM by applying the true daily withdrawals to the respective replenishment plans. This allows for a counterfactual comparison, evaluating the proposed strategy by assuming that the demand that arose in reality would have occurred even if the strategy had been operationally implemented. 

The performance of both the QUBO model and the benchmark solution is measured using two primary Key Performance Indicators (KPIs): total replenishment cost and average service level.

\begin{enumerate}
    \item \textbf{Total Replenishment Cost.} The first KPI measures the total operational cost associated with all replenishment activities scheduled by the model. The cost of a single replenishment for a given machine $i$ includes both fixed and variable components. The total cost for the entire benchmark period is the sum of the individual costs for every replenishment scheduled by the model.
    \item \textbf{Average Service Level.} The second KPI, Average Service Level (ASL), measures the ability of the machines to meet customer withdrawal demand, i.e., the fraction of demand that is actually satisfied. The most correct metric would be to consider the actual time in which the machine remains without cash, but being in a simulation environment, this information is not available and for this reason we must consider an estimate.
 For each ATM $i$ on day $t$, let $\hat{R}_{i,t}$ be the cash available at the start of the day based on the simulation strategy and let $W^{*}_{i,t}$ be the real withdrawal demand observed that day. The unmet demand (lost withdrawals due to cash shortage) is defined as 
    \begin{equation}{\text{lost}}_{i,t}=\max\{0,\,W^{*}_{i,t}-\hat{R}_{i,t}\}
    \end{equation}
    In fact if $\hat{R}_{i,t}>W^{*}_{i,t}$ all the withdrawal $W^{*}_{i,t}$ would be dispensed; otherwise, if $\hat{R}_{i,t}\le W^{*}_{i,t}$, all the residual $\hat{R}_{i,t}$ would be dealt.
    Aggregating over all ATMs and all days, the average service level is defined as 
    \begin{equation}
    ASL=1-\sum_{i,t}\frac{\text{lost}_{i,t}}{W^{*}_{i,t}}
    =1-\sum_
    {i,t}\frac{\max\{0,\,W^{*}_{i,t}-\hat{R}_{i,t}\}}{W^{*}_{i,t}}
    \end{equation}
    This definition yields $ASL=1$ when all demand is satisfied (no cash-out events) and decreases as the amount of unmet demand increases.
\end{enumerate}

\subsection{Results}\label{subsec44}
The QUBO-based optimization model was benchmarked against the threshold-based solution using the KPIs defined in Section \ref{subsec43}. The comparison was conducted over four representative months: April, May, October, and November 2022. The results demonstrate significant advantages of the optimization approach in terms of cost efficiency, while maintaining a high standard of service.

\begin{table}[h]
\caption{Service level KPIs: Threshold-based vs. QUBO Optimization Strategy}\label{tab:lds}%
\begin{tabular}{@{}lcccccccc@{}}
\toprule
\textbf{KPI: Service level} & \textbf{Solution} & \textbf{Apr} & \textbf{May} & \textbf{Oct} & \textbf{Nov} \\
\midrule
\textbf{Average} & Threshold-based & 99.91\% & 99.90\% & 99.71\% & 99.74\% \\
& QUBO & 99.92\% & 99.93\% & 99.78\% & 99.82\% \\
\midrule
\textbf{Minimum} & Threshold-based & 99.42\% & 99.37\% & 99.26\% & 98.17\% \\
& QUBO & 99.00\% & 99.47\% & 96.99\% & 98.44\% \\
\botrule
\end{tabular}
\end{table}

\begin{table}[h]
\caption{Cost KPIs: Threshold-based vs. QUBO Optimization Strategy}\label{tab:cost}%
\begin{tabular}{@{}lcccccccc@{}}
\toprule
\textbf{KPI: Cost in \euro} & \textbf{Solution} & \textbf{Apr} & \textbf{May} & \textbf{Oct} & \textbf{Nov} \\
\midrule
\textbf{Total Recharging} & Threshold-based & 60\,943 & 60\,224 & 79\,343 & 76\,871  \\
& QUBO  & 50\,384 & 50\,986 & 64\,868 & 64\,575 \\
\midrule
\textbf{Average per Day} & Threshold-based & 2\,031  & 1\,943  & 2\,559  & 2\,562\ \\
& QUBO & 1\,680  & 1\,645  & 2\,093  & 2\,153 \\
\midrule
\textbf{Saving} & & \textbf{17.33\%} & \textbf{15.34\%} & \textbf{18.24\%} & \textbf{16.00\%} \\
\botrule
\end{tabular}
\end{table}

Service level KPIs are reported in Table~\ref{tab:lds}. The optimization solution consistently achieves an~\emph{Average Service Level} that is either comparable to or slightly higher than the threshold-based solution across all simulated months. This indicates that the model's primary goal of cost reduction does not come at the expense of overall machine availability. For instance, in May, the service level improved from 99.90\% to 99.93\%, and in November it rose from 99.74\% to 99.82\%.
The~\emph{Minimum Service Level}, which reflects the worst-case performance for any single ATM, shows a more complex trade-off. Although the optimization model improved this metric in May and November, it registered a lower minimum level in April and notably in October (96.99\% vs. 99.26\%). This suggests that while the model effectively manages the network on average, its aggressive cost-saving strategy may, in specific scenarios, allow individual machines to approach their safety stock levels more closely than the conservative threshold-based method. This highlights a classic trade-off between optimality and risk, which can be tuned via the model's penalty parameters.

Cost KPIs are reported in Table~\ref{tab:cost}. The most significant outcome of the benchmark is the substantial reduction in operational costs achieved by the QUBO model. The optimization solution delivered consistent savings across all four months, ranging from \textbf{$15.3\%$} to \textbf{$18.2\%$}. In absolute terms, this translates to monthly savings between approximately \euro $\,9\,200$ and \euro $\,14\,500$ for $276$ ATMs. For example, in October, total costs were reduced from \euro $\,79\,343$ to \euro $\,64\,868$, an improvement of $18.24\%$. Similarly, the average daily cost of operations was markedly lower in all periods, demonstrating the model's ability to generate more efficient and economically advantageous replenishment schedules.

In conclusion, the experimental results validate the effectiveness of the QUBO framework for cash flow optimization. The model consistently provides significant double-digit cost savings while maintaining an excellent average service level, proving its superiority over the threshold-based methodology.
\section{Conclusions}\label{sec5}
This paper has introduced a QUBO-based approach to optimizing cash replenishment for networks of Automated Teller Machines, addressing the joint challenges of demand uncertainty, nonlinear cost structures, and stringent service level requirements. The key modeling elements include an objective function that encodes fixed, variable, and co-location discount components, feasibility penalties enforcing at most one replenishment per ATM in the planning horizon and a service level penalty that preserves compatibility with the QUBO structure through precomputed shortage costs.

On the methodological side, we show that the intrinsic structure of ATM cash logistics—binary “replenish or not” decisions over a rolling multi-day horizon, quadratic interactions due to co-location discounts, and constraints expressible via penalties—maps naturally to the QUBO framework. The resulting QUBO instances are solved using MegaQUBO, a Data Reply's proprietary GPU-based QUBO solver, which exploits massive parallelism on modern GPU architectures to handle realistically sized networks and planning horizons within practical computation times. This combination of a problem-specific QUBO formulation with an industrial-grade GPU solver turns the proposed methodology into a practically deployable decision-support tool for financial institutions, rather than a purely theoretical proof of concept.

The empirical evaluation, based on a dataset of 276 ATMs in Italy over four representative months of 2022, demonstrates the effectiveness of the proposed approach. Compared with the threshold-based policy, this solution achieves consistent double-digit reductions in total recharging costs—between about 15\% and 18\%—while maintaining an average service level close to threshold-based solution. These results confirm that a QUBO formulation, when solved on a high-performance GPU-based solver, can substantially improve cost efficiency without materially compromising customer service quality, although some trade-offs are observed in minimum service levels for individual ATMs when the model is tuned aggressively for cost savings.


This study has already been translated and used into a production-grade decision-support pipeline, integrating rolling-horizon optimization with operational data feeds and an upstream forecasting component, and confirming cost reduction under real execution frictions (e.g., delays and partial service). 
Beyond direct financial savings, the proposed approach has broader operational and ESG implications. In fact, reducing the number and frequency of replenishment operations limits vehicle use and associated emissions, decreases cash handling in processing centers, and can contribute to safer and leaner cash operations.

Given the favorable computational behavior of the approach, we are evaluating extensions to a larger operational perimeter that includes cash recycling ATMs. These machines introduce an additional layer of complexity, as the cash deposited by customers is recirculated into the dispensing pool, thereby requiring the forecasting module to produce two independent demand estimates—one for withdrawals and one for incoming deposits—whose net aggregation determines the effective daily cash flow. Quantifying how cost savings and service level KPIs evolve as the network size increases and the forecasting task becomes inherently more challenging constitutes a natural direction for future work.
\bmhead{Acknowledgements}

The authors would like to thank the domain Gestione Integrata Dei Valori and Experimentations and Labs of Innovation and Business Support Rengineering at Intesa Sanpaolo, as well as all those involved in the project whose contributions helped turn this proof of concept into a live application.

\bibliography{sn-bibliography}


\end{document}